\documentclass[10pt]{article}
\usepackage[dvips]{graphicx}
\usepackage{amsmath}
\usepackage{amsthm}
\usepackage{pstricks}
\usepackage{amssymb}
\usepackage{mathtools}
\usepackage{mathrsfs}
\usepackage[all]{xy}
\usepackage{hyperref}
\usepackage{dsfont}

\theoremstyle{plain}
\newtheorem{theorem}{Theorem}

\newtheorem{remark}{Remark}

\theoremstyle{remark}

\theoremstyle{definition}

\newcommand{\Ham}[0]{{\rm Ham}}
\newcommand{\Sp}[0]{{\rm Symp}}
\newcommand{\Diff}[0]{{\rm Diff}}
\newcommand{\ham}[0]{{\mathfrak{ham}}}
\newcommand{\symp}[0]{{\mathfrak{sp}}}

\title{Symplectic parallel integrators in the realm of Hofer's geometry}
\author{Hugo Jim\'enez-P\'erez}
\date{\today}
\begin{document}
\maketitle
\begin{abstract}
  Symplectic integrators constructed from Hamiltonian and Lie formalisms
  are obtained as symplectic (indeed Hamiltonian) maps whose flow follows
  the exact solution of a ``sourrounded'' Hamiltonian $\tilde H=H+h^kH_1$. 
  Those modified Hamiltonians depends virtually on the time by the timestep 
  size $h$.  When the numerical integration of a Hamiltonian 
  system involves more than one symplectic scheme as in the 
  \emph{parallel-in-time} algorithms and specifically the Parareal scheme,
  there are not a simple way to control the dynamical behavior of the error
  Hamiltonian. The interplay of to different symplectic integrators can
  degenerate their behavior if both have different dynamical properties, 
  reflected in the number of iterations to have a good approximation to 
  the final sequential solution. Considered as flows of time-dependent 
  Hamiltonians we use the Hofer's geometry to search for the optimal 
  coupling of symplectic schemes.  As a result,
  we obtain the constraints in the Parareal method to have a good behavior 
  for Hamiltonian dynamics.
\end{abstract}

\section{Introduction}
Symplectic integrators are the natural methods for simulating
Hamiltonian dynamics. 
To construct a symplectic integrator, we can follow two different procedures:
on one side we can use the Hamiltonian formalism using generating functions,
Lie transforms, etc. On the other side we can take a general method, for 
instance, Runge-Kutta formulas, modifying the coefficients to satisfy the 
symplecticity conditions \cite{Abd1}. In general, we obtain implicit methods but when the Hamiltonian 
can be separated into kinetic and potential energies $H(q,p)=T(p)-V(q)$ we can 
construct explicit methods easy to implement.

They are constructed using the 
diffeomorphisms which left invariant the symplectic form 
$\omega = dp\wedge dq$ which defines
de Hamiltonian vector field $X_H$.
Such a diffeomorphisms are called symplectic diffeomorphisms or 
symplectomorphisms and they form a subgroup denoted by $\Sp(M,\omega)$.
In particular, the flow of any Hamiltonian vector field is a symplectic 
diffeomorphism and all of them form another subgroup of diffeomorphisms 
called the Hamiltonian diffeomorphisms, denoted by $\Ham(M,\omega)$.
%We have the contentions
%\begin{eqnarray*}
%  \Ham(M,\omega)\subset \Sp(M,\omega)\subset \Diff(M).
%\end{eqnarray*}

Moreover, using the Lie formalism we can consider $\Sp(M\omega)$ and $\Ham(M,\omega)$
as Lie groups and the sets of symplectic and Hamiltonian vector fields as
their Lie algebras at the identity element. Then, for every Hamiltonian vector
field $X_H\in\mathfrak{ham}(M,\omega)$, its flow is given by the exponential 
$\varphi_t(x_0)= e^{tX_H}x_0$. For a fixed $t=h\in\mathbb R^+$ the mapping 
$x(t_0)\mapsto x(t_0+h)$ generated by $\varphi_h=e^{hX_H}$ is a Hamiltonian 
map which defines the symplectic integrator. In fact, this is a Hamiltonian 
integrator which preserves more structure than the symplectic one.

Classically the paralelization of this type of systems is
performed by a decomposition of the phase space (domain decomposition) or 
looking for parallelizable tasks into the method or into the equations.
However, in the last two decades there were several attempts to develop 
another type of parallelization for symplectic integrators based on the 
decomposition of the time variable. The first attempt of some 
parallel-in-time algorithm for a scalar differential equation 
was published by Nievergelt in 1964 \cite{Nie1}. The idea is to 
decompose the total time in several subintervals 
which can be computed in parallel. Each interval,
called a ``branch'', must be modified propagating the corrected local initial
condition to each subinterval; this technique has derived in the 
\emph{multishooting} methods. Although
the Nievergelt's algorithm is not iterative, almost all other algorithms
use an iterative process to approximate the sequential numerical solution.
Those iterative algorithms consist in two steps: one \emph{predictor} which 
estimates in parallel the value of several branches and one 
\emph{corrector}\footnote{Many authors inverse the terminology 
considering the parallel step as the corrector and the sequential step 
as the predictor. Our choise is evident when we relate the corrector step
with the \emph{symplectic correctors} studied in \cite{WHT96,McL1,MP00}} 
which approaches the final solution propagating the 
predictions between different branches.
We call them the \emph{time-parallel} algorithms and all of them
differ in the corrector step which uses different iterative process 
to convergence. Of course, there are others differences but in this paper we 
are interested in the corrector step. 

As noted by Saha, Stadel and Tremaine \cite{SST97} one way to time-parallelize 
an almost integrable Hamiltonian system is to compute in parallel
several branches saving the perturbing contributions and to propagate
them computing in sequence the integrable Hamiltonian part. This technique 
is reproduced in \cite{JL11} for high-order symplectic integrators however,
it works fine if the ratio of the computing time of the integrable over
the perturbing part is very small. This approach is very accurate but
expensive. On the other extreme, there is the \emph{parareal}
method introduced by Lions, Maday and Turinici in \cite{LMT01}
and refined in \cite{BM02}. In this method the propagation of the 
predictions is made by a simple increment of the corrector at every 
iteration. As a result, it is a very fast algorithm, however, for 
 Hamiltonian systems there are several inconvenients 
associated to the non preservation of the geometric structure 
of the underlying integrators. 

In order to deal with this type of problems, Bal and Wu \cite{BW08} have done 
the first step considering a new way to spread the information between 
branches in the sequential step and practically destroying the
``pure'' parareal scheme. Also, Dai \emph{et al.} \cite{DLLM10}
have introduced another variation of the parareal step, using symmetries 
and projections into the energy manifold to preserve the geometric 
properties of the underlying integrators. Recently, the author has 
proposed a geometric corrector step using Lie's algebras \cite{Jim2}
which is equivalent to that from Dai \emph{et. al.}. 
In this paper we translate that approach to the Hofer's geometry in
order to search not only for Hamiltonian maps but for optimizing 
the energy (which we relate with the number of iterations) to go from 
the first guess solution to the final solution. 

\section{Lie algebras and Hamiltonian vector fields}
We consider the phase space of a Hamiltonian system 
as a symplectic manifold\footnote{All the computations and 
results listed here apply to arbitrary symplectic manifolds
$(M,\omega)$.}  $M=T^*\mathbb R^n\cong\mathbb R^{2n}$
with the canonical symplectic form $\omega=dp\wedge dq$.
Denote by $\mathfrak X(M)$ the set of all vector fields
and by $\mathfrak F(M)\cong C^\infty(M)$ the
set of all smooth functions over $M$.

We define de binary operation $[\cdot,\cdot]:\mathfrak X(M)
\times \mathfrak X(M)\to\mathfrak X(M)$ by
the rule
\begin{eqnarray}
  \left[ X, Y \right] &=& XY -YX,\qquad X,Y\in\mathfrak X(M),
  \label{eqn:liebrack:def}
\end{eqnarray}
called the Lie bracket, which is: bilinear, alternating, and satisfies
the Jacobi identity.
The set of all the vector fields $\mathfrak X(M)$ 
equipped with the Lie bracket (\ref{eqn:liebrack:def}) obtains the 
structure of Lie algebra $(\mathfrak X(M),[\cdot,\cdot])$.

Let $\mathcal L_XF$  be the Lie derivative of $F$ along the vector 
field $X\in\mathfrak X(M)$.
$\mathcal L_XF$ meassures the change of $F$ along 
$X$ where $F$ can be a function, a vector field, 
a $p$-form or, in general, a tensor.

The Lie derivative of the symplectic form $\omega$ along $X$ is given by the Cartan's 
magic formula
\begin{eqnarray}
  \mathcal L_X\omega = di_X\omega + i_Xd\omega,
  \label{eqn:Lie:der}
\end{eqnarray}
where $d$ is the exterior differential and $i_X\omega=\omega(X,\cdot)$
is the contraction of $\omega$ by $X$ or equivalently the inner 
product of the vector field $X$ with the 2-form $\omega$.

We say that the vector field $X$ is \emph{symplectic} if its flow preserves
$\omega$, which means
$\mathcal L_X\omega=0$.
Since 
$\omega$ is the canonical symplectic 2-form then $\omega = -d\lambda$ where 
$\lambda=pdq$ is the Liouville form. Consequently, the second term in the
right hand side of (\ref{eqn:Lie:der}) is zero.
In other words, a vector field $X$ is symplectic if the 1-form $i_X\omega$ 
is closed which means that it belongs to the kernel of $d$
\begin{eqnarray*}
  d\left( i_X\omega \right) = d\left( \omega(X,\cdot) \right)=0.
\end{eqnarray*}
We denote the set of symplectic vector fields on $M$ by 
$\symp(M,\omega)$.

We say that $X$ is \emph{Hamiltonian} if, 
in addition, $i_X\omega$ is exact, \emph{i. e.}, there exists 
$f\in\mathfrak F(\mathbb R^{2n})$ such that
\begin{eqnarray}
  i_X\omega = \omega(X,\cdot) = -df.
  \label{eqn:Ham:def}
\end{eqnarray}
We call $f$ a Hamiltonian function for $X$ and we write $X=X_f$ to 
specify that $X$ is the Hamiltonian vector field associated to $f$. 
the set of all Hamiltonian vector fields on $M$ is denoted by $\ham
(M,\omega)$. Finally, the triplet $(M,\omega,X_H)$ defines a Hamiltonian 
system over $M$. For general mechanical systems the configuration 
space can be consider as a Riemannian manifold $(N,g)$ and the 
phase space becomes the cotangent bundle $M=T^*N$ which has
a natural structure of symplectic manifold.

We need a non obvious result from the theory of differential 
$p$-forms which give us the expression of the inner product 
of the Lie bracket with a $p$-form $\alpha$ over any differential 
manifold $M$.
The inner product $i_{[X,Y]}\alpha$ for any $p$-form $\alpha$ is 
given by (\cite[pp 73]{Ber1}):
\begin{eqnarray*}
  i_{[X,Y]}\alpha = \mathcal L_X(i_Y\alpha) - i_Y(\mathcal L_X\alpha).
\end{eqnarray*}

Consequently, for every two $X,Y\in\symp(M,\omega)$ we have
\begin{eqnarray*}
  i_{[X,Y]}\omega &=&  \mathcal L_X(i_Y\omega)+0\\
   &=& d\left( i_X\left( i_Y\omega \right) \right)+0\\
   &=& d(\omega\left( Y,X \right))\\
   &=& -d\left( \omega\left( X,Y \right) \right)
\end{eqnarray*}
where we used $\mathcal L_Y\omega =0$ and $d\omega=0$.

These computations has important consequences: 1) $[X,Y]$ is a Hamiltonian 
vector field. 2) Since $\ham(M,\omega)\subset\symp(M,\omega)$, the Lie bracket 
$[\cdot,\cdot]$ gives them the structure of Lie subalgebras of 
$\mathfrak X(M)$, and
3) $\ham(M,\omega)$ is an (in fact the maximal) ideal 
of $\symp(M,\omega)$ with respect to $[\cdot,\cdot]$, \emph{i.e.}
\begin{eqnarray*}
  [\symp(M,\omega),\symp(M,\omega)]\subset \ham(M,\omega).
\end{eqnarray*}
%It means that for every $X,Y\in\mathfrak{sp}(\mathbb R^{2n})$ then 
%$[X,Y]\in\mathfrak{ham}(\mathbb R^{2n})$. 
We have the relations
\begin{eqnarray}
  \ham(M,\omega)\subset \symp(M,\omega)\subset \mathfrak X(M).
   \label{eqn:seq:alg}
\end{eqnarray}
It is a well-known fact that $\ham(M,\omega) = \symp(M,\omega)$ if and only 
if the fundamental group of $M$ is trivial; in other words, when $M$ is 
simply connected.

Now we link this point of view with the classical development in local
coordinates. Let $H: M\to\mathbb R$ be a differentiable function
with Hamiltonian vector field $X_H$ on $M$. Select a point 
$v\in \phi^{-1}(M)$ in a local chart of M. The Darboux's theorem 
says that, locally, all symplectic manifolds are symplectomorphic
to $T^*\mathbb R^n\cong \mathbb R^{2n}$ and then we can consider 
that $v\in\mathbb R^{2n}$. In canonical symplectic 
coordinates $v=(q,p)$ such that 
$q\in \mathbb R^n$ and $p\in T^*_q\mathbb R^n\cong\mathbb R^n$. 
The vector field $X_H$ in local coordinates is
\begin{eqnarray}
  \dot q = \frac{\partial H}{\partial p},&\qquad&
  \dot p = -\frac{\partial H}{\partial q},
  \label{eqn:pq}
\end{eqnarray}
which are called the Hamilton equations. 

For two differentiable functions $f,g\in\mathfrak F(M)$, their 
associated Hamiltonian vector fields $X_f, X_g$ fullfils 
$\omega\left( X_f,X_g \right) =-\{f,g\}$ where $\{\cdot,\cdot\}$ 
is the \emph{Poisson bracket} for functions defined by
\begin{eqnarray}
  \left\{ f,g \right\}&:=& 
     \frac{\partial f}{\partial q} \frac{\partial g}{\partial p} 
     - \frac{\partial g}{\partial q} \frac{\partial f}{\partial p}.
     \label{eqn:pb:def}
\end{eqnarray}
The binary operation (\ref{eqn:pb:def}) is bilinear, antisymmetric
and fulfills the Jacobi identity. 
The space of real-valued differentiable functions $\mathfrak F(M)$,
equipped with the Poisson bracket  (\ref{eqn:pb:def}), obtains the 
structure of Lie algebra. It is possible to write the Hamiltonian 
vector field in terms of the Poisson brackets by
\begin{eqnarray*}
  \dot z = X_H(z) = \left\{ z, H \right\}.
\end{eqnarray*}

There exists a natural anti-morphism of Lie algebras between 
the algebra of differentiable functions $\mathfrak F(M,
\{\cdot,\cdot\})$ and the algebra of Hamiltonian vector fields 
$\mathfrak X(M,[\cdot, \cdot])$ given by
\begin{eqnarray}
    f&\mapsto& X_f\\
    \{f,g\}&\mapsto& X_{\{f,g\}}= - \left[ X_f, X_g \right] 
  \label{eqn:morph}
\end{eqnarray}
Since the Poisson bracket of two functions is a function
then the Lie bracket of two Hamiltonian vector fields 
is again a Hamiltonian vector field, as we have shown before.

\section{Hamiltonian diffeomorphisms and Hofer's geometry}

%One way to study the general properties of mathematical and 
%geometrical objects is by the study of its symmetries and invariants.
%Differential geometry is the study of the transformations which 
%preserves the differential structure between geometrical objects. 
%In the same way, symplectic and Hamiltonian geometries study the 
%geometry of objects which are invariants under symplectic 
%and Hamiltonian diffeomorphims, respectively.

A symplectic diffeomorphism or symplectomorphism of a symplectic 
manifold $(M,\omega)$ is a $C^\infty$ diffeomorphism $\phi\in\Diff(M)$ 
which preserves the symplectic structure $\omega$, it means that the 
\emph{pull-back} of $\phi$ fulfills $\phi^*(\omega)=\omega$.
The support of a diffeomorphism $\phi$ is the closure of 
$\{x\in M |\phi(x)\neq x\}$. 
The set of all symplectomorphisms with compact support form a group, 
denoted $\Sp(M, \omega)$  (with the law of composition of mappings).

Let $X_H\in\ham(M,\omega)$ be a Hamiltonian vector field. The flow 
$\varphi^H_t=e^{tX_H}(x_0)$ of $X_H$ is a one parameter 
subgroup of symplectic diffeomorphisms\footnote{We use the 
exponential map of vector fields since $\Sp(M,\omega)$ is actually
a Lie group}. The set of all the symplectomorphisms which arise
as the flow $\varphi^f_t$ of Hamiltonian vector fields $X_f$ form 
another subgroup of diffeomorphisms called the Hamiltonian 
diffeomorphisms, denoted by $\Ham(M,\omega)$. It is easy to prove that
$\Ham(M,\omega)$ is an infinite-dimensional subgroup since for every
$f\in\mathfrak F(M)$, we have a mapping $X_\_:\mathfrak F(M)\to
\mathfrak X(M)$ which maps $f\mapsto X_f$ with 
${\rm ker}(X_\_)\cong\mathbb R$ the constant functions. 
The exponential map $e^\_:\ham(M,\omega)\to\Ham(M,\omega)$ 
is injective then the composition $e^\_\circ X_\_(H)=e^{tX_H}$ 
sends an infinite
dimensional basis of $\mathfrak F(M)$ to an infinite-dimensional 
basis of $\Ham(M,\omega)$.

We can relate the groups $\Ham$, $\Sp$ and $\Diff$ with the Lie algebras
in (\ref{eqn:seq:alg}) by the exponential map as follows
\begin{eqnarray}
    \begin{array}{ccccc}
  \Ham(M,\omega) &\subset & \Sp(M,\omega)& \subset& \Diff(M)\\
     \exp \uparrow &      & \exp \uparrow &      & \exp \uparrow\\
  \ham(M,\omega) &\subset& \symp(M,\omega)&\subset& \mathfrak X(M).
  \end{array}
   \label{eqn:seq:grp}
\end{eqnarray}
The reader must note that the exponential mapping is, in general,
not surjective and the group $\Sp(M,\omega)$ can have several
components. By construction $\Ham(M,\omega)$ belongs to the identity 
component of $\Sp(M,\omega)$.

Consider the set of time-dependent Hamiltonian functions $H:M\times I
\to\mathbb R$ with compact support. We can \emph{normalize} such a 
functions since for every interval $I=[0,a], a<\infty$ the flow of
the vector field associated to the function $\tilde H(x,t)=aH(x,at)$ 
defined on $\tilde I=[0,1]$ is again a Hamiltonian flow. More 
generally for every smooth function $g(t)$ with $g(0)=0$ the flow
$\varphi^H_{g(t)}$ is Hamiltonian with normalized Hamiltonian
function $\frac{d g}{dt}(t)H(x,g(t))$ \cite{Pol1}. These properties 
of rescaling  in time are used to regularize singularities in 
mechanical systems.  

We can restate the definition of Hamiltonian diffeomorphism in the 
following way:  a \emph{Hamiltonian diffeomorphism} is a 
diffeomorphism $\phi:M\to M$ which can be written as the 
time-1-map of a Hamiltonian flow, \emph{i.e.}, $\phi = \varphi^H_1$
for some time-periodic Hamiltonian $H : \mathbb S^1 \times M\to R$
\cite{Sib1}. Let us denote by $H_t$ the function $H(t,\cdot)$ on 
$M$. In the following, we will normalize the Hamiltonians and consider 
only the time-1-maps in $\Ham(M,\omega)$. We denote by $\mathcal 
H(M,\omega)$ the group of time-1-map of Hamiltonian flows which is
a subgroup of $\Ham(M,\omega)$.

\begin{remark}
   There is a constraint in all these definitions since the theory
applies for functions and diffeomorphisms with compact support. 
However, for symplectic integrators, we do not need global 
properties and the most important thing is the numerical tests
for the error behavior.
\end{remark}

A diffeomorphism $\phi\in\Diff(M)$ is said to be \emph{isotopic 
to the identity} if there exists a smooth
map $H : M \times [0, 1]\to M$ such that if $h_t: M\to M$ is given 
by $h_t(x) = H(x, t)$, then $h_t$ is a $C^\infty$ diffeomorphism,
$h_0 = id_M$ and $h_1 = \phi$. We say that $h_t$ is an isotopy
from $\phi$ to the identity.

We say that $h_t$ is a \emph{Hamiltonian isotopy} if there exists a 
smooth family of functions $H_t : M\to \mathbb R$ such that
\begin{eqnarray}
    i_{X_{ht}}\omega = -dH_t
\end{eqnarray}

In this way, we have constructed curves or trayectories in 
$\Ham(M,\omega)$ which connects any Hamiltonian diffeomorfism 
in $\mathcal H(M,\omega)$ with the identity map.

For every $\phi\in \Ham(M, \omega)$, choose a Hamiltonian isotopy 
$\Phi = (\phi_t)$ from $\phi$ to the identity.       
Hofer \cite{HZ94} defined the length of this isotopy by
\begin{eqnarray}
 l_H(\Phi) :=\int_{\mathbb S^1} {\rm osc}( H_t) dt 
\end{eqnarray}
where ${\rm osc}(H_t) := \max(H_t) - \min(H_t)$ denotes the oscillation of a 
function on M. For $H \in \mathcal H(M,\omega)$, it is clear 
that $l(H) = 0$ if, and only if, $H = 0$. 
The \emph{distance from the identity}, or \emph{energy}, of an 
element $\phi\in\Ham(M, \omega)$ is defined as
\begin{eqnarray}
    d(id, \phi) := \inf\{  l(H) | H \in \mathcal H, \phi = 
   \varphi^H_1 \}.
   \label{eqn:met:hof}
\end{eqnarray}
Let us extend the distance to a function 
$d : \Ham(M, \omega)\times\Ham(M, \omega)\to[0, \infty)$
by setting $d(\phi, \psi) := d(id, \psi\circ \phi ^{-1} )$.
Hofer has showed in \cite{Hof1} that $d(\cdot,\cdot)$ is a 
bi-invariant metric on $\Ham(M,\omega)$ defined intrinsecally.
Then the set $\Ham(M,\omega)$ with the metric (\ref{eqn:met:hof})
is called the Hofer's geometry and it is a fundamental stone in 
symplectic topology.

As was pointed out by Siburg in \cite{Sib1}, to every Hamiltonian
dynamical system corresponds one single path in $\Ham(M,\omega)$ and
vice versa. All the dynamical properties of the Hamiltonian system,
including its periodic orbits, heteroclinic connections, etc. are 
contained in the Hamiltonian isotopy. Moreover, Baily and Polterivich
have showed in \cite{BP1} that the bifurcation diagram of 
every Hamiltonian system is preserved for every Hamiltonian diffeomorphism
contained in a geodesic path in $\Ham(M,\omega)$. 

For instance, let $H_t \in\mathcal H(M,\omega)$ be an admissible Hamiltonian
which means that $H_t$ has a compact support for every $t\in I$. Then $H_t$
is said to \emph{generate a minimal geodesic} if $d(id,\varphi_t^H ) = l(H)$.
However, it is very difficult to work with the space of geodesics in
$\Ham(M,\omega)$.

\section{Symplectic integrators and Hamiltonian maps}

Symplectic integrators are the natural numerical methods for simulating
Hamiltonian dynamics. From the geometrical point of view, the most natural
integrators are obtained by the Hamiltonian formalism and the Lie theory
applied to the group $\Ham(M,\omega)$ and its Lie algebra 
$\ham(M,\omega)\cong T_{Id}\left(\Ham(M,\omega)\right)$.

Consider the Hamiltonian system $(M,\omega,X_H)$. The flow generated
by the Hamiltonian vector field $X_H\in\ham(M,\omega)$ is the one-parameter
subgroup of $\Ham(M,\omega)$ defined by $\varphi^H_t(x_0)=e^{tX_H}(x_0)$
where $x_0\in M$ is the initial condition of the vector field.
For a fixed value $t=\tau$ the mapping $e^{\tau X_H}:M\to M$ is a symplectic, 
actually a Hamiltonian map. In the generic case, $e^{\tau X_H}$ is very 
complicated and is given in (analytical) implicit form. 

However, there exists an important class of Hamiltonian functions for which 
the Legendre condition\footnote{The Legendre condition ask for the convexity 
of the Hamiltonian function for which its Hessian does not vanish.} is 
satisfied. This class is formed by separable functions $H(q,p)=T(p)+V(q)$
where $T(p)$ is the kinetic energy and $V(q)$ is the potential. (We use
the plus sign in order to develop the exponential as a product.) Then 
the Hamiltonian vector field is separable and each part can be integrated
independently obtaining an implicit method. If we write $X_H=X_T+X_V$
the flow generated by $X_H$ becomes
\begin{eqnarray*}
    \varphi_t^H = e^{\tau X_H} = e^{\tau (X_T+X_V)}
\end{eqnarray*}
and since each part can be integrated independently, we can estimate
$e^{\tau X_T}$ and $e^{\tau X_V}$ directly. However, in general $[X_T,X_V]\neq 
0$ and therefore $$e^{\tau(X_T+X_V)}\sim e^{\tau X_T}e^{\tau X_V}$$ 
concide only in the first term. This becomes a first order method.

In order to obtain higher order methods, we search for coefficients 
$\{a_i\}_1^m$ and $\{b_i\}_1^m$ such that $\sum_i a_i = 1$ and $\sum_i b_i=1$
to estimate $e^{\tau X_H}$ by the composition of maps
\begin{eqnarray}
    e^{\tau X_T + \tau X_V} \sim e^{a_1\tau X_T} e^{b_1\tau X_V}\dots 
          e^{a_m\tau X_T} e^{a_m\tau X_V}.
	  \label{eqn:hord}
\end{eqnarray}
We must impose additional conditions to the coefficients $a_i$ and $b_i$
for matching more terms on both sides of the expression (\ref{eqn:hord}). 
For that, we use the Baker-Campbell-Hausdorff (BCH) formula and derive a set 
of polynomial conditions for $a_i$ and $b_i$. In particular, the derivation of 
time-parallel methods impose a reversibility condition which requires that 
the symplectic method be symmetric \cite{Jim2}.

Using the BCH formula we can write the symplectic integrator (\ref{eqn:hord})
of order $k$ with its residue as
\begin{eqnarray}
    e^{\tau X_T + \tau X_V} = e^{a_1\tau X_T} e^{b_1\tau X_V}\dots 
    e^{a_m\tau X_T} e^{a_m\tau X_V} + \mathcal O(\tau^{k+1}).
	  \label{eqn:hord:ex}
\end{eqnarray}
which implies 
\begin{eqnarray}
    e^{a_1\tau X_T} e^{b_1\tau X_V}\dots 
    e^{a_m\tau X_T} e^{a_m\tau X_V} = e^{\tau X_H} - \mathcal O(\tau^{k+1}).
	  %\label{eqn:hord:ex}
\end{eqnarray}
The numerical solution follows a modified or ``surrounded'' 
Hamiltonian which we consider as the nonautonomous function
\begin{eqnarray}
   \tilde H(t,q,p)=H(q,p) + t^{k+1}H_1(t,q,p). 
  \label{eqn:sur:Ham}
\end{eqnarray}
This property give us a measure for the error commited by the symplectic
methods which  corresponds to the error Hamiltonian $\tau^{k+1}H_1(\tau,q,p)$ 
for fixed $\tau$. The interested readers will find a deeper discussion on 
the subject in the 
works of Suzuki \cite{Suz1}, Yoshida \cite{Yos1} for this type of 
Hamiltonians and McLachlan \cite{McL1}, Laskar and Robutel \cite{LR01} 
for Hamiltonians of the form $H(q,p)=T(p)+ \epsilon V(q)$ for small 
$\epsilon$. A more general discussion on the construction of symplectic 
integrators using the BCH formula is found in Hairer \emph{et.al.}
\cite{HLW1}.

For every symplectic integrator obtained with this process, we have a 
Hamiltonian function (\ref{eqn:sur:Ham}) characterized by its 
error Hamiltonian which depends on $(\{a_i\},\{b_i\},\tau)$, and then 
an analytical flow $\varphi^{\tilde H}\in\Ham(M,\omega)$ associated to 
the numerical solution. For 
$f,g\in\mathfrak F_0(M)$ such that $f\neq g$ the Hamiltonian maps
\begin{eqnarray}
  \varphi^f = e^{\tau X_f}\quad{\rm and}\quad \varphi^g = e^{\tau X_g},
    \qquad \varphi^f,\varphi^g\in\mathcal H(M,\omega),
  \label{eqn:flows}
\end{eqnarray}
associated to their flows fulfills $\varphi^f\neq \varphi^g$.
The Hofer's metric give us a way to measure the distance between them
in $\Ham(M,\omega)$. This fact permit us to compare symplectic 
integrators in an geometrical framework. 
%In fact, we are interested in the 
%displacement energy and geodesics between Hamiltonian maps.

\section{Time-parallel integrators}
The corner stone in the theory of time-parallel integrators is the 
use of two different numerical flows, generally using a two level 
discretization $\delta t \prec \Delta t \prec [0,T]$ and a corrector step
which propagates sequentially with the coarse flow, the values 
obtained in parallel with the fine flow. 
In the seminal article of 
Nievergelt \cite{Nie1} he proposed, as an example, using the Euler 
method on each discretization. This produces two numerical flows 
which interplays to estimate the final sequential solution.

The first documented proposal for a time-parallel algorithm for 
Hamiltonian systems (using a symplectic integrator), is in the work 
of Saha, Stadel and Tremain \cite{SST97}. They propose a time-parallel
algorithm for almost integrable Hamiltonian systems in action-angle
coordinates. In such type of coordinates, the integrable part of 
the system corresponds to the actions which are constants 
in the flow. Instead of having a two
level discretization, they use a single level discretization 
with two flows: one for the integrable system and the other one
for the complete (perturbed) system using the symplectic mid-point 
rule. In this work, at least from the examples they show, it is not 
possible to use a high order symplectic method.

The next interesting proposal is the Parareal algorithm which was 
used for Hamiltonian systems with a lot of degrees of freedom
\cite{AMV09, CF06, FC03}. However, for long time simulations of 
systems with a few degrees of freedom Parareal has not a 
well behavior \cite{FC03,FHM05}. The problem is that the parareal
algorithm (in fact, the parareal step) does not preserves the 
symplecticity when it propagates the fine flow with the coarse one.
There are several attempts to obtain a better corrector step
in order to preserve the symplecticity \cite{BW08, DLLM10,Jim2}, 
however, there are not a concrete answer to this problem.
Here we ask for an additional point.

Given two different symplectic
maps, one for each discretization we can obtain the
symplectic map of its Lie bracket and construct a symplectic 
corrector as in \cite{BW08,Jim2}. At each iteration corresponds 
a point in $\Ham(M,\omega)$ and then there exists a Hamiltonian
isotopy $h_t\subset\Ham(M,\omega)$  which contains all such points. 
Since the fine flow is given (it corresponds to the final numerical
solution we expect to approximate), the question rests in the 
choice of the coarse flow in order to approach the fine solution
in the minimal number of iterations. We give a partial answer in the 
next paragraphs.

Let's consider the initial value problem
\begin{eqnarray}
  \dot y(t)= X_H(y(t)),\qquad y(0)=y_0.
  \label{eqn:sys:1}
\end{eqnarray}
where $y:[0,T]\to M$ and $X_H\in\ham(M,\omega)$.
We discretize the problem 
by partitioning $[0,T]$ in $N$ subintervals of size $\Delta t = T/N$
which we call \emph{branches} as in \cite{Nie1}.
We set $t_0 = 0$, $t_n = T$ and $t_i=i\Delta t$ such that 

\begin{eqnarray}
   0=t_0<t_1<t_2< \cdots < t_i <\cdots <t_N=T.
\end{eqnarray}
We write $\Delta t\prec [0,T]$ for this partition. 
Each branch in $\Delta t$ is decomposed in $N_\delta\in\mathbb N$ 
subintervals $\delta t=\Delta t/N_\delta$
which corresponds to the final resolution. Then 
$\delta t\prec\Delta t\prec [0,T]$.

On this discretization we introduce two levels of resolution: a 
coarse symplectic solver $G$ on the $\Delta t$ resolution and a 
fine symplectic solver $F$ on the $\delta t$ one.
Their flows $\varphi^G$ and $\varphi^F$ are defined 
uniquely by their coefficients and their timesteps
$(\{a_i\},\{b_i\}, \Delta t)$ and $(\{c_i\},\{d_i\}, \delta t)$
respectivelly.

The fine and the coarse propagation of 
the solutions $y_n$ on a branch are respectively 
given by $\mathcal F(y_n)$ and $\mathcal G(y_n)$. 
For the problem (\ref{eqn:sys:1}) we estimate 
the first guess sequence with the coarse solver $y_0^{(0)}=y(0)$
and $y_{n+1}^{(0)}= \mathcal G(y_{n}^{(0)})$. The time-parallel
algorithms are given by some variation of the to steps
\begin{eqnarray}
  y_{0}^{(k+1)} &=& y_0^{(k)},\\
  y_{n+1}^{(k+1)}&=&  
           \Gamma\left( \mathcal F\left(y_{n}^{(k)}\right),
            \mathcal G\left(y_{n}^{(k+1)}\right),
            \mathcal G\left(y_n^{(k)}\right)\right),
  \label{eqn:sys1}
\end{eqnarray}
where the subscripts are the propagation in time and the superscripts
are the iterations. In (\ref{eqn:sys1}), $\Gamma$ is the corrector 
which is an operator on the symplectic solvers.
The general algorithm is as follows
\begin{center}
  \begin{tabular}[ht]{|rl|}
    \hline
    \multicolumn{2}{|l|}{
    {\bf Time-parallel algorithm 1.}
    }\\
    \hline
    1: & Setup of the initial guess sequence\\
    2: & $\quad$ $y_0^0=y(0)$, $y_{n+1}^0=\mathcal G\left(y_n^0\right)$\\
    3: & For $k=1$ to $k_{max}$\\
    4: & $\quad$ \emph{Parallel} resolution on $[T^n,T^{n+1}]$: \\ 
    5: & $\qquad$ compute $\mathcal F \left(y_{n}^{(k)}\right)$,\\
    6: & $\quad$ For $n=k$ to N\\
    7: & $\qquad$ \emph{Sequential} corrections: \\
    8:   & $\quad\qquad$ compute $y_{n+1}^{(k+1)}= 
           \Gamma\left(\mathcal F\left(y_{n}^{(k)}\right),
            \mathcal G\left(y_{n}^{(k+1)}\right),
            \mathcal G\left(y_n^{(k)}\right)\right)$.\\
    9: & $\quad$end for (n).\\
    10: & end for (k).\\
    \hline
  \end{tabular}
\end{center}

In particular the parareal algorithm implements (\ref{eqn:sys1}) as
\begin{eqnarray}
  y_{n+1}^{(k+1)}&=&  
            \mathcal F\left(y_{n}^{(k)}\right)
           + \mathcal G\left(y_{n}^{(k+1)}\right)
           - \mathcal G\left(y_n^{(k)}\right).
  \label{eqn:sys11}
\end{eqnarray}
Expression (\ref{eqn:sys11}) has been called \emph{the parareal 
iteration}.

Recall that the construction of a symplectic corrector 
uses the Lie bracket of $X_F$ and $X_G$ with at least one reversal 
integration, then we impose that both symplectic solvers be 
symmetrical \cite{Jim2}.

Every implicit symplectic integrator for separable Hamiltonian 
systems can be uniquely determined by its coefficients $\{c_i\}$
and $\{d_i\}$. Its numerical flow $\varphi^{\tilde H}$ by the
coefficients and the timestep $(\{c_i\},\{d_i\},\tau)$. Then, we 
identify the flow of every symplectic integrator with timestep
$\tau$ with the triplet $\varphi^{\tilde H}\cong(\{c_i\},\{d_i\},\tau)$.
Suppose that we have, as the fine symplectic scheme of order $k$, 
the more accurated integrator
in the family of the $k$-order symplectic integrators.  

\begin{theorem}
    Given the two level discretization $\delta t\prec\Delta t\prec 
    [0,T]$ and the final (optimal) symplectic integrator 
    $(\{c_i\},\{d_i\},
    \delta t)$ with flow $\varphi^F$ the closest flow $\varphi^G$
    in $\Ham(M,\omega)$ for $\Delta t$ is given by $(\{c_i\},
    \{d_i\},\Delta t)$.

    \label{teo:1}
\end{theorem}
{\it Proof.} Immediate using the triangle's inequality property 
of the Hofer's metric.
$\hfill\square$

\begin{figure}[ht]
    \begin{center}
	\includegraphics[scale=0.65]{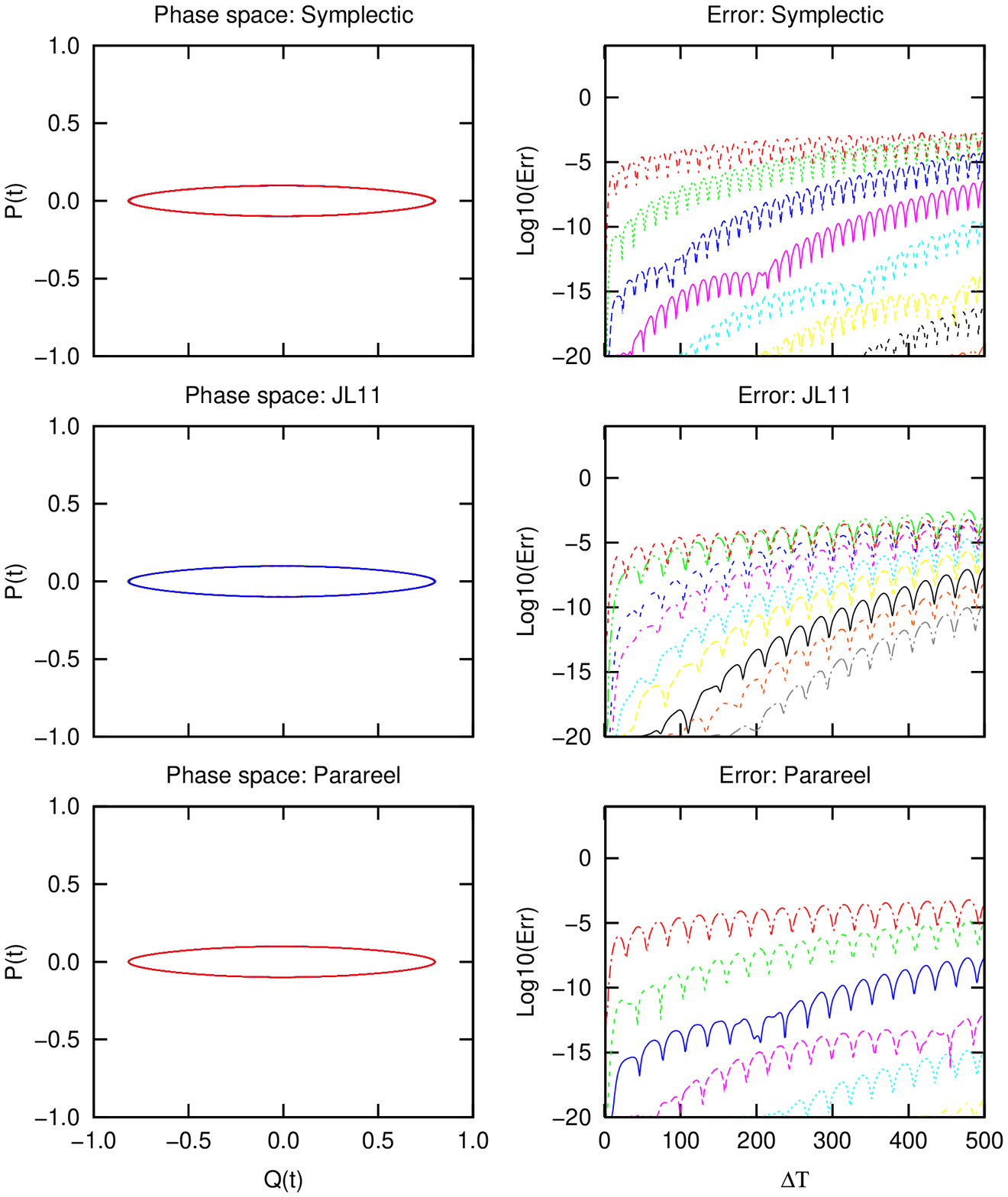}
    \end{center}
    \caption{Numerical test of the Spin-orbit problem with JL11 and 
    parareal. The parammeters are $\varepsilon=0.1$, $\alpha=0.01$,
    $\theta=0.2$, $\delta t=1/128$, $\Delta t=1$, $p_0=0$, $q_0=0.8$. }
    \label{fig:1}
\end{figure}

As a consequence, the use of some symplified symplectic scheme 
increases the number of iterations and, in general, introduces 
an erratic behavior since the dynamics of both flows is in  general
not equivalent.

In this case, the corrector step constructed by the author in
\cite{Jim2} coincides with 
those studied by Wisdom and Holman \cite{WHT96} McLachlan \cite{McL2}
and Laskar and Robutel \cite{LR01}.

There are more results from the symplectic topology which apply for 
particular families of Hamiltonians. For example, a function $H\in
\mathfrak F_0(M,\omega)$ is called to have \emph{quadratic growth}
if there exists $c\in\mathbb R$ finite such that  $|d^2H(z)|<c$
for all $z\in M$. For such functions, the Hamiltonian isotopy 
of the Theorem \ref{teo:1} is a geodesic in $\Ham(M,\omega)$ 
\cite[Ch. 12]{MS1}.

\section{A numerical test}
We compare the pure parareal with the more accurate JL11 algorithm
introduced in \cite{JL11} for the 
one dimensional Spin-orbit problem with several values of the 
parameter $\alpha$  
\begin{eqnarray}
    H(q,p) &=& \frac{1}{2}p^2 -\varepsilon
    \left( \cos(2q) + \alpha \left( \cos 2q+\theta - 7 \cos 2q 
    -\theta \right) \right)
    \label{eqn:Ham}
\end{eqnarray}
We select a 8th order symplectic integrator  from the 
$\mathcal{SBAB}_n$ family for both $\mathcal F$ and $\mathcal G$ 
solvers. 

The first approximation with the coarse solver approach faster
than the symplectic proposed in \cite{Jim2} and JL11 proposed in 
\cite{JL11}. However, as was showed in the later document, JL11
obtain exactly the same sequential solution.

\end{document}